\documentclass[12pt,a4paper]{article}
\usepackage{amssymb,amsmath}
\usepackage{amsfonts}
\usepackage{a4}
\begin{document}
 \baselineskip18pt
 \title{\bf  Good integers : A note on results of  Jitman and  Prugsapitak}
\author{ Madhu Raka
 \\ \small{\em Centre for Advanced Study in Mathematics}\\
\small{\em Panjab University, Chandigarh-160014, INDIA}\\
\date{}}
\maketitle
 {\abstract{In this paper we fix some errors made by  Jitman \cite{Jit} and Prugsapitak and Jitman \cite{PJit} while characterizing good integers and $2^{\beta}$-good integers.}\vspace{2mm}
 }
 \section{ Introduction.} For fixed coprime nonzero integers $a$ and $b$, a  positive integer $\ell$ is called a good integer
with respect to $a$ and $b$ (see  Moree \cite{Mor}), if there exists a positive integer $k$ such that $\ell|(a^k + b^k)$. Otherwise, $\ell$ is called a bad integer. Denote
by $G_{(a,b)}$  the set of good  integers defined with
respect to $a$ and $b$. A positive integer $\ell$ is said to be oddly-good (with respect to $a$ and $b$) if
$\ell|(a^k + b^k)$ for some odd integer $k \geq 1$, and evenly-good  if $\ell|(a^k + b^k)$
for some even integer $k \geq 2$. Therefore, $\ell$ is good if it is oddly-good or evenly-good. Denote
by $OG_{(a,b)}$ (resp., $EG_{(a,b)}$) the set of oddly-good (resp., evenly-good) integers. For a non-negative
integer $\beta$, a positive integer $d$ is said to be $2^{\beta}$-good (with respect to $a$ and $b$) if $2^{\beta}d \in G_{(a,b)}$.
Otherwise, $d$ is said to be $2^{\beta}$-bad. In the same fashion,
$2^{\beta}$-oddly-good  and $2^{\beta}$-evenly-good integers are defined.  For an integer $\beta\geq 0$, denote by $G_{(a,b)}(\beta), OG_{(a,b)}(\beta)$ and $EG_{(a,b)}(\beta)$ the
sets of $2^{\beta}$-good, $2^{\beta}$-oddly-good, and $2^{\beta}$-evenly-good integers, respectively.\vspace{2mm}

In \cite{Jit}, Jitman characterized good integers, oddly-good integers and considered their applications in coding theory. But there are some  errors in this paper.  In the proof of Proposition 2.1 of \cite{Jit}, it is used that for   odd integers $a$ and $b$,
\begin{equation} Ord_{2^{\beta}}(ab^{-1}) =2 ~~~ \Rightarrow ~~~ ab^{-1} \equiv -1 ( {\rm mod}~ 2^{\beta}) ~~~~  {\rm i. e.} ~~~2^{\beta}\mid a+b.\end{equation}
Again in the proof of Proposition 2.3 of \cite{Jit}, it is used that
\begin{equation} Ord_{d}(ab^{-1}) =2k ~~~ \Rightarrow ~~~ (ab^{-1})^k \equiv -1 ( {\rm mod}~ d),\end{equation}
where $a$ and $b$ are coprime to $\ell=2^{\beta}d$, $\beta \geq 1$, $d$ is an odd positive integer and $b^{-1}$ denotes the multiplicative inverse of $b$ modulo $\ell$.\vspace{2mm}

\noindent These are false statements as $$ Ord_8(11)=2, ~~ 11^2 \equiv 1 ( {\rm mod}~ 8) ~~  {\rm but }~~ 11 \not\equiv -1( {\rm mod}~8).$$
$$ Ord_{15}(11)=2, ~~ 11^2 \equiv 1 ( {\rm mod}~ 15) ~~  {\rm but } ~~ 11 \not\equiv -1( {\rm mod}~15).$$

 Because of these errors, Proposition 2.1, Proposition 2.3, Theorem 2.1, Corollary 2.1, Theorem 3.1 and Corollary 3.3 of \cite{Jit} are no longer true. Proofs of some of otherwise correct results also need to be modified. \vspace{2mm}

 In a subsequent paper \cite{PJit}, the authors Prugsapitak and Jitman tried to fix the second error (though not mentioning the error explicitly), but they still overlooked the first error. Because of this, Proposition 2.1, Proposition 2.2, Proposition 2.3, Corollary 2.1 and Corollary 2.2  of \cite{PJit}  are again no longer true.  As a consequence, applications of good integers in the study of self-dual negacyclic codes are also affected. In fact the statement that $$ Ord_{d}(ab^{-1}) =2k ~~~ \Rightarrow ~~~ (ab^{-1})^k \equiv -1 ( {\rm mod}~ d)$$
 holds only when $d$ is an odd prime power or $d=2$. The aim of this paper is to fix these errors and to rectify the above mentioned propositions and results.
 \section{ Rectified results.}
Through out the paper let $a$, $b$ and $\ell=2^{\beta}d $, where $\beta \ge 0$ and $d$ an odd positive integer, be pairwise coprime non-zero integers.  Let $b^{-1}$ denote the multiplicative inverse of $b$ modulo $\ell$ and $Ord_m(ab^{-1})$ denote the multiplicative order of $ab^{-1}$ modulo $m$ for  a divisor $m$ of $\ell$. It is clear that $\ell \in G_{(a,b)}$ or $ d \in G_{(a,b)}(\beta)$ if and only if $(ab^{-1})^k \equiv -1 ( {\rm mod}~ \ell)$ for some positive integer $k$. Let $ x= ab^{-1}$. Denote by $2^{\gamma}||\ell$ if $\gamma $ is the largest integer such that $2^{\gamma}|\ell$. Note that if gcd$(a,b)=1$ and $\ell \in G_{(a,b)}$ then gcd$(a,\ell)=1$ and gcd$(b,\ell)=1$. \vspace{2mm}


\noindent {\bf Lemma 1}: Let $p$ be an odd prime and $r$ be a positive integer. Let $p^r$ be  good and $s$ be the smallest positive integer such that $(ab^{-1})^s \equiv -1 ( {\rm mod}~ p^r)$.  Then $Ord_{p^r}(ab^{-1})=2s$. \vspace{2mm}

\noindent This is  Proposition 2 of \cite{Mor}.\vspace{2mm}

\noindent The converse of Lemma 1 is also true. Let $ x= ab^{-1}$. If $Ord_{p^r}(x)=2s$, we have $p^r | (x^s-1)(x^s+1)$. It can not happen that $p^{i} | (x^s-1)$ and $p^j | (x^s+1)$ with $i+j=r, ~i \geq 1, j\ge 1$. Because then $p| (x^s-1)$ and $p | (x^s+1)$ which gives $p | 2$; not possible as $p$ is an odd prime.  Hence  either
$p^r | (x^s-1)$ or $p^r | (x^s+1)$ but not both. If  $p^r | (x^s-1)$, we get $Ord_{p^r}(x) \geq s$, not possible. Therefore $p^r$ must divide $(x^s+1)$. \vspace{2mm}

\noindent In fact, we have a more general result.\vspace{2mm}

\noindent {\bf Lemma 2}: Let $a$, $b$ and $d$  be pairwise coprime odd integers. If $k$ is the smallest positive integer such that $(ab^{-1})^k \equiv -1 ( {\rm mod}~ d)$  then $Ord_{d}(ab^{-1})=2k.$ \vspace{2mm}

\noindent {\bf Proof}: Let $k= 2^{\lambda} k'$, $\lambda \ge 0, ~k'$ odd. Let $x= ab^{-1}$.   As $x^k \equiv -1 ( {\rm mod}~ d)$ we have $x^{2k} \equiv 1 ( {\rm mod}~ d)$. Therefore $Ord_{d}(x)\mid 2k.$ Let $Ord_{d}(x)= r = 2^{\mu} r'$, where  $0\leq \mu \le \lambda+1, ~r'$ is odd and $ r'| k'$. Let $k'=r'r''$.\vspace{2mm}

\noindent If $\mu \leq \lambda $, $x^{2^{\mu}r'} \equiv 1 ( {\rm mod}~ d)$, gives $x^{2^{\lambda}k'} \equiv (x^{2^{\mu}r'})^{2^{\lambda-\mu}r''} \equiv 1 ( {\rm mod}~ d)$, but $x^{2^{\lambda}k'} \equiv x^k \equiv -1 ( {\rm mod}~ d)$. Therefore $1\equiv -1 ( {\rm mod}~ d)$. This is not possible as $d$ is odd. Therefore we must have $\mu =\lambda + 1$. \vspace{2mm}

\noindent If $d >1$, let $d = p_1^{e_1}p_2^{e_2}\cdots p_t^{e_t}$ where $p_i$ are odd primes and $e_i \geq 1$. As $x^r = x^{2^{\lambda +1}r'} \equiv 1 ( {\rm mod}~ p_i^{e_i})$ we have $p_i^{e_i} | (x^{2^{\lambda}r'}-1)(x^{2^{\lambda}r'}+1)$ for each $i$. As before, it can not happen that $p_i^{\alpha_i}| (x^{2^{\lambda}r'}-1)$ and $p_i^{\beta_i} | (x^{2^{\lambda}r'}+1)$ for some $\alpha_i \geq 1, \beta_i \ge 1$ with $\alpha_i+\beta_i=e_i$.   Hence for each $i,~ 1 \le i \le t$,  either
$p_i^{e_i} |(x^{2^{\lambda}r'}-1)$ or $p_i^{e_i} | (x^{2^{\lambda}r'}+1)$ but not both. If for some $i$, $p_i^{e_i} | (x^{2^{\lambda}r'}-1)$, we get $x^k  \equiv  x^{2^{\lambda}k'} \equiv x^{2^{\lambda}r'r''}  \equiv  1 ( {\rm mod}~ p_i^{e_i})$. Not possible as we are given that $ x^k  \equiv -1 ( {\rm mod}~ p_i^{e_i})$ and $p_i$ is odd. \vspace{2mm}

\noindent Hence $p_i^{e_i} \mid (x^{2^{\lambda}r'}+1)$ for all $i$. Therefore $d \mid (x^{2^{\lambda}r'}+1)$ i.e. $x^{r/2} \equiv -1 ( {\rm mod}~ d)$. Now the minimality of $k$ gives $r/2=k$.  ~~~~~~~~~~~~~~~~~~~~~~~~~~~~~~~~~~~~~~~~~~~~~~~~~~~~~~~$\Box$\\

\noindent  The converse of Lemma 2 is not always true as illustrated in Section 1. Note that
$$ Ord_{2^{\beta}}(x)=\left \{ \begin{array}{ll} 1 & {\rm if} ~ \beta=1\\ 2 & {\rm if} ~ \beta \geq 2 ~{\rm and }~ x\equiv -1 ( {\rm mod}~ 2^{\beta}). \end{array}\right.$$
If $ \ell= 2^{\beta}p_1^{e_1}p_2^{e_2}\cdots p_t^{e_t}$ where $p_i$ are odd primes and $\beta \geq 0,~ e_i \geq 0$, we have
$$ Ord_{\ell}(x)= {\rm lcm}\big( Ord_ {2^{\beta}}(x), Ord_ {p_1^{e_1}}(x), Ord_ {p_2^{e_2}}(x), \cdots, Ord_ {p_t^{e_t}}(x) \big).$$

\noindent Following are some results of Moree \cite{Mor}.

\noindent {\bf Lemma 3} (\cite{Mor}, Proposition 2): For an odd prime $p$, $Ord_{p^{e}}(x)= Ord_p(x)p^{\alpha}$ for some $\alpha\geq 0$. \vspace{2mm}

\noindent {\bf Lemma 4} (\cite{Mor}, Theorem 1): Let $d>1$ be an odd integer. Then $ d \in G_{(a,b)}$ if and only if there exists an integer $s\geq 1$ such that $2^s|| Ord_p(x)$ for every prime $p$ dividing $d$. \vspace{2mm}

\noindent {\bf Lemma 5} (\cite{Jit}, Proposition 2.2): Let $a,b,d>1$ be pairwise coprime odd integers. Then $ d \in G_{(a,b)}$ if and only if $ 2d \in G_{(a,b)}$. \vspace{2mm}

\noindent The correct form of Proposition 2.1 of \cite{Jit} and Proposition 2.2 of \cite{PJit} is \vspace{2mm}

\noindent {\bf Proposition 1}: If $a$, $b$ are coprime odd integers and $\beta\geq 1$ is any integer, then the following are equivalent :\vspace{2mm}

$ \begin{array}{ll} (1)& 2^{\beta}|a+b\\(2)&2^{\beta}\in G_{(a,b)}\\ (3)&\beta=1 ~{\rm or~} ab^{-1} \equiv -1 ( {\rm mod}~ 2^{\beta}). \end{array}$\vspace{2mm}

\noindent The correct form of Proposition 2.3 of \cite{Jit} and Proposition 2.1 and Corollary 2.1 of \cite{PJit} is \vspace{2mm}

\noindent {\bf Proposition 2}: Let $a,b,d>1$ be pairwise coprime odd integers and $\beta \geq 2$ be any integer. Then $2^{\beta}d \in G_{(a,b)}$ if and only if $ab^{-1} \equiv -1 ( {\rm mod}~ 2^{\beta})$ and $2|| Ord_p(ab^{-1})$ for every prime $p$ dividing $d$.\vspace{2mm}

\noindent{\bf Proof}: Suppose $2^{\beta}d \in G_{(a,b)}$. Let $k$ be the smallest positive integer such that $2^{\beta}d|(a^k + b^k)$. This gives
$2^{\beta}|(a^k + b^k)$. If $k$ is even $$ a^k+b^k = (a^2)^{k/2}+(b^2)^{k/2}\equiv 1+1 \equiv 2 ( {\rm mod}~ 4),$$ as an odd square is always congruent $1$ modulo $4$. Therefore $k$ must be odd. But then $$ a^k+b^k =(a+b)\big(a^{k-1}-a^{k-2}b+a^{k-3}b^2- \cdots + b^{k-1}\big).$$ The second factor on the right hand side is odd, it being a sum of odd terms taken odd number of times. Therefore $2^{\beta}|(a + b)$ which gives $ab^{-1} \equiv -1 ( {\rm mod}~ 2^{\beta})$. Also $k$ is smallest integer such that $d|(a^k + b^k)$, i. e., $x^k \equiv -1 ( {\rm mod}~ d)$. Then we have, by Lemma 2, $Ord_{d}(x)=2k,$ where $k$ is odd. Let $d = p_1^{e_1}p_2^{e_2}\cdots p_t^{e_t}$ where $p_i$ are odd primes and $e_i \geq 1$. Then, using Lemma 3,
$$ 2k=Ord_d(x)= {\rm lcm}\big(  Ord_{p_1}(x)p_1^{\alpha_1}, Ord_{p_2}(x) p_2^{\alpha_2}, \cdots, Ord_{p_t}(x)p_t^{\alpha_t} \big).$$
Also $x^k \equiv -1 ( {\rm mod}~ p_i)$ for all $i, 1\leq i \leq t$ and $k$ is odd. Therefore $Ord_{p_i}(x)$ is even and $2|| Ord_{p_i}(x)$ for each $i$.

Conversely let $2|| Ord_{p_i}(x)$ for each $p_i| d$. This gives $2|| Ord_{p_i^{e_i}}(x)$ for each $i$. Let $ Ord_{p_i^{e_i}}(x)=2r_i$, where $r_i$ is odd. Therefore $x^{r_i} \equiv -1 ( {\rm mod}~ p_i^{e_i})$ for all $i, 1\leq i \leq t$. Let $k= {\rm lcm}(r_1,r_2.\cdots,r_t)$, $k$ is odd and let $k=r_ir'_i$. Each of $r'_i$ is also odd. Then $x^k \equiv x^{r_i r'_i} \equiv (-1)^{r'_i} \equiv -1( {\rm mod}~ p_i^{e_i})$ for each $i$. Therefore $x^k \equiv -1 ( {\rm mod}~ d)$. Now $x \equiv -1 ( {\rm mod}~ 2^{\beta})$ implies $x^k \equiv -1 ( {\rm mod}~ 2^{\beta})$ as $k$ is odd. Hence $x^k \equiv -1 ( {\rm mod}~ 2^{\beta}d)$, i.e., $2^{\beta}d \in G_{(a,b)}$.\vspace{2mm}

In view of the above results, Theorem 2.1, Corollary 2.1, Theorem 3.1 and Corollary 3.3 of \cite{Jit} should  read as follows :

\noindent {\bf Theorem 1}: Let $a$ and $b$ be pairwise coprime non-zero integers and let $\ell=2^{\beta}d$ be a positive integer such that $d$ is odd and $\beta \geq 0$.
\begin{enumerate} \item If $ab$ is odd, then $\ell = 2^{\beta}d \in G_{(a,b)}$ if and only if one of the following statements hold \vspace{2mm}\\ $\begin{array}{ll} {\rm (a)}
 &  \beta \in \{0,1\} {\rm ~and ~} d=1\vspace{2mm}\\ {\rm (b)} &   \beta \in \{0,1\},~ d \geq 3 {\rm ~and ~there~ exists} ~ s\geq 1 {\rm ~ such~ that~} 2^s|| Ord_p(ab^{-1}) {\rm ~ for}  \\ & {\rm ~every~ prime~} p {\rm ~ dividing~} d.\vspace{2mm}\\{\rm(c)} &  \beta \geq 2, d=1 {\rm ~ and~} ab^{-1} \equiv -1 ( {\rm mod}~ 2^{\beta}).\vspace{2mm}\\ {\rm(d)}&  \beta \geq 2,~ d\geq 3,~ ab^{-1} \equiv -1 ( {\rm mod}~ 2^{\beta}) {\rm ~ and~} 2|| Ord_p(ab^{-1}) {\rm ~ for ~every~ prime} \\ & p {\rm ~ dividing~} d.\end{array}$
\item If $ab$ is even, then $\ell = 2^{\beta}d \in G_{(a,b)}$ if and only if one of the following statements hold \vspace{2mm}\\ $ \begin{array}{ll} {\rm (a)}& \beta =0 {\rm ~ and ~} d=1.\vspace{2mm}\\{\rm (b)}& \beta =0, ~ d \geq 3 {\rm ~ and ~there~ exists~} s\geq 1 {\rm ~ such ~ that~} 2^s|| Ord_p(ab^{-1}) {\rm ~ for ~every} \\ & {\rm ~prime~} p {\rm ~ dividing~} d. \end{array} $\end{enumerate}

\noindent {\bf Corollary 1}: Let $a,~b$ and $\ell$ be pairwise coprime non-zero integers and let $\ell=2^{\beta}d$ be a positive integer such that $d$ is odd and $\beta \geq 0$. Let $\gamma \geq 0 $ be an integer such that $2^{\gamma}||a+b$. Then $\ell  \in G_{(a,b)}$  if and only if one of the following statements hold.
\begin{enumerate} \item $\ell =1,2$. \item $d=1$ and $2\le \beta \leq \gamma$.\item
$d\geq 3$, $\beta \in \{0,1\}$ and $2^s|| Ord_p(ab^{-1})$  for some $s \geq 1$ and for every prime $p$ dividing $d$. \item $d\geq 3$, $2\leq \beta \leq \gamma$ and $2|| Ord_p(ab^{-1})$ for every prime $p$ dividing $d$.
\end{enumerate}
In 3 and 4, if $\ell  \in G_{(a,b)}$ then $2^s|| Ord_{\ell}(ab^{-1})$ if and only if $2^s|| Ord_p(ab^{-1})$ for every prime $p$ dividing $d$.\vspace{2mm}

\noindent {\bf Theorem 2}: Let $a$ and $b$ be pairwise coprime non-zero integers and let $\ell=2^{\beta}d$ be a positive integer such that $d$ is odd and $\beta \geq 0$.
\begin{enumerate} \item If $ab$ is odd, then $\ell = 2^{\beta}d \in OG_{(a,b)}$ if and only if one of the following statements hold \vspace{2mm}\\ $\begin{array}{ll} {\rm (a)}
 &  \beta \in \{0,1\} {\rm ~and ~} d=1\vspace{2mm}\\ {\rm (b)} &   \beta \in \{0,1\},~ d \geq 3 {\rm ~and ~} 2|| Ord_p(ab^{-1}) {\rm ~ for}  {\rm ~every~ prime~} p {\rm ~ dividing~} d.\vspace{2mm}\\{\rm(c)} &  \beta \geq 2, d=1 {\rm ~ and~} ab^{-1} \equiv -1 ( {\rm mod}~ 2^{\beta}).\vspace{2mm}\\ {\rm(d)}&  \beta \geq 2,~ d\geq 3,~ ab^{-1} \equiv -1 ( {\rm mod}~ 2^{\beta}) {\rm ~ and~} 2|| Ord_d(ab^{-1}) \end{array}$
\item If $ab$ is even, then $\ell = 2^{\beta}d \in OG_{(a,b)}$ if and only if one of the following statements hold \vspace{2mm}\\ $ \begin{array}{ll} {\rm (a)}& \beta =0 {\rm ~ and ~} d=1.\vspace{2mm}\\{\rm (b)}& \beta =0, ~ d \geq 3 {\rm ~ and ~}  2|| Ord_p(ab^{-1}) {\rm ~ for ~every}  {\rm ~prime~} p {\rm ~ dividing~} d. \end{array} $\end{enumerate}

    \noindent {\bf Corollary 2}: Let $a,~b$ and $\ell$ be pairwise coprime non-zero integers and let $\ell=2^{\beta}d$ be a positive integer such that $d$ is odd and $\beta \geq 0$. Let $\gamma \geq 0 $ be an integer such that $2^{\gamma}||a+b$. Then $\ell  \in OG_{(a,b)}$  if and only if  one of the following statements hold.
\begin{enumerate} \item $\ell =1,2$. \item $d=1$ and $2\le \beta \leq \gamma$.\item
$d\geq 3$, $0\leq \beta \le \gamma$ and $2|| Ord_p(ab^{-1})$   for every prime $p$ dividing $d$. In that case  $2|| Ord_{\ell}(ab^{-1})$.
\end{enumerate}

\end{document}